%% file: pMFinal.tex
\documentclass[12pt]{amsart}
\usepackage{amsmath, amssymb, latexsym, amsthm, mathrsfs, color, mathtools}
\usepackage[hyphens]{url}
\usepackage[T1]{fontenc}
\usepackage{lmodern}

\usepackage[top=2in, bottom=1.5in, left=1in, right=1in]{geometry}

\input{shortcuts}

\nc{\ed}{

\subsection*{Acknowledgments}
We are indebted to Diego Alejandro Mej\'\i{}a for his Lemma~\ref{lem:mej}.
We also thank Will Brian and Ashutosh Kumar for their answers to questions
we had during this study~\cite{Bri, Ash}, and Heike Mildenberger
for information about the hypothesis $\fr<\fd$. 
The research of the first named author was supported by \emph{Etiuda~2},
Polish National Science Center, UMO-2014/12/T/ST1/00627.

\end{document}
}

\title[Products of Menger spaces]{Products of Menger spaces:\\
a combinatorial approach}

\author[P. Szewczak]{Piotr Szewczak}
\address{Piotr Szewczak, 
Institute of Mathematics, Faculty of Mathematics and Natural Science College of Sciences, Cardinal Stefan Wyszy\'nski University in Warsaw, W\'oycickiego 1$\slash$3, 01--938 Warsaw, Poland,
and
Department of Mathematics, Bar-Ilan University, Ramat Gan 5290002, Israel
}
\email{p.szewczak@wp.pl}

\author[B. Tsaban]{Boaz Tsaban}
\address{Boaz Tsaban,
Department of Mathematics, Bar-Ilan University, Ramat Gan 5290002, Israel}
\email{tsaban@math.biu.ac.il}
\urladdr{http://math.biu.ac.il/~tsaban}

\begin{document}

\begin{abstract}
We construct Menger subsets of the real line whose product is not 
Menger in the plane. In contrast to earlier constructions,
our approach is purely combinatorial.
The set theoretic hypothesis used in our construction is far milder than earlier
ones, and holds in all but the most exotic models of real set theory.
On the other hand, we establish productive properties
for versions of Menger's property parameterized by filters and semifilters.
In particular, the Continuum Hypothesis
implies that every productively Menger
set of real numbers is productively Hurewicz, and
each ultrafilter version of Menger's property is strictly
between Menger's and Hurewicz's classic properties.
We include a number of open problems emerging from this study.
\end{abstract}



\maketitle

\section{Introduction}

A topological space $X$ is \emph{Menger}
if for each sequence
$\eseq{\cU}$
of open covers of the space $X$, there are finite subsets
$\cF_1\sub\cU_1$, $\cF_2\sub\cU_2$, \dots
whose union forms a cover of the space $X$.
This property was introduced by Karl Menger~\cite{Menger24}, and reformulated as
presented here by Witold Hurewicz~\cite{Hure25}.
Menger's property is strictly between $\sigma$-compact and Lindel\"of.
Now a central notion in topology, it has applications in a number of branches of
topology and set theory. The undefined notions in the following example,
which are available in the indicated references,
are not needed for the remainder of this paper.

\bexm
Menger spaces form the most general class for which
a positive solution of the D-space problem is known~\cite[Corolarry~2.7]{AurD}, and
the most general class for which
a general form of Hindman's Finite Sums Theorem holds~\cite{AlgSelRT}.
In set theory, Menger's property characterizes filters whose Mathias forcing
notion does not add dominating functions~\cite{ChoZdo}.
\eexm

Menger's property is hereditary for closed subsets and continuous images.
By a classic result of Todor\v{c}evi\'c there are,
provably, Menger spaces $X$ and $Y$ such that the
product space $X\x Y$ is not Menger~\cite[\S3]{Tod95}.
It remains open whether there are, provably, such examples in the real line, or even just metrizable examples~\cite[Problem~6.7]{OPiT}.
This problem, proposed by Scheepers long ago, resisted tremendous efforts thus far.

For brevity, sets of real numbers are called here \emph{real sets}.
Assuming \CH{}, there are two Luzin sets whose product is not Menger~\cite[Theorem~3.7]{coc2}.
An uncountable real set is \emph{Luzin} if its intersection with every
meager (Baire first category) set is countable. An uncountable real set $X$ is \emph{concentrated} if
it has a countable subset $D$ such that the set $X\sm U$ is countable for
every open set $U$ containing $D$. Every Luzin set is concentrated, and
every concentrated set has Menger's property.
This approach extends to obtain similar examples using a set theoretic hypothesis about the meager sets that is weaker than \CH{}~\cite[Theorem~49]{cbc}.
Later methods~\cite[Theorem~9.1]{sfh} were combined with reasoning
on meager sets to obtain examples using another portion of 
\CH{}~\cite[Theorem~3.3]{RZ}.
Here, we introduce a purely combinatorial approach to this problem. We obtain examples
using hypotheses milder than earlier ones, as well as examples using hypotheses
that are incompatible with \CH{}. To this end, we introduce 
the key notion of \emph{bi-$\fd$-unbounded set}, and determine 
the limits on its possible existence.
We extend these results to variations of Menger's property parameterized by
filters and semifilters.
For a semifilter $S$, we introduce the notion of \emph{$S$-scale}.
These scales provably exist, and capture a number of
distinct special cases used in earlier works.

The second part of the paper, beginning with Section~\ref{sec:pR}, establishes
provably productive properties among semifilter-parameterized Menger properties.
If $S$ is an \emph{ultrafilter}, then every $S$-scale gives rise to a productively $S$-Menger space.
We deduce that each of these variations of Menger's property is strictly between Hurewicz's
and Menger's classic properties.

The last section includes a discussion of related results and open problems suggested
by our results.

\section{Products of Menger spaces}

Towards a combinatorial treatment of the questions discussed here, we identify
the Cantor space $\Cantor$ with the family $\PN$ of all subsets of the set $\bbN$.
Since the Cantor space is homeomorphic to Cantor's real set, every
subspace of the space $\PN$ is considered as a real set.

The space $\PN$ splits into two important subspaces: the family of infinite
subsets of $\bbN$, denoted $\roth$, and the family of finite subsets of $\bbN$, denoted $\Fin$.
We identify every set $a\in\roth$ with its increasing enumeration, an element of the Baire
space $\NN$.
Thus, for a natural number $n$, $a(n)$ is the $n$-th element in the increasing
enumeration of the set $a$.
This way, we have $\roth\sub\NN$, and the topology of the space $\roth$ (a subspace
of the Cantor space $\PN$) coincides
with the subspace topology induced by $\NN$. This explains some
of the elementary assertions made here; moreover,
notions defined here for $\roth$ are often adaptations of classic notions for $\NN$.
Depending on the interpretation,
points of the space $\roth$ are referred to as sets or functions.

For functions $a,b\in\roth$, we write $a\le b$ if $a(n)\le b(n)$ for all natural numbers
$n$, and
$a\les b$ if $a(n)\le b(n)$ for almost all natural numbers $n$, that is,
the set of exceptions $\sset{n}{b(n)<a(n)}$ is finite.
We follow the convention that \emph{bounded}
means \emph{has an upper bound} in the ambient superset.

\bdfn
Let $\kappa$ be an infinite cardinal number.
A set $X\sub\roth$ with $\card{X}\ge\kappa$ is \emph{$\kappa$-unbounded}
if the cardinality of every $\le$-bounded subset of the set $X$ is smaller than $\kappa$.
\edfn

\brem
For cardinal numbers $\kappa$ of uncountable cofinality, which will be the case in the present paper,
the notion of $\kappa$-unbounded defined here is equivalent to its variation using the relation
$\le^*$ instead of $\le$.
This is not the case for cardinal numbers of countable cofinality.
\erem

Let $\kappa$ be an infinite cardinal number. A topological space $X$ with $\card{X}\ge\kappa$ 
is \emph{$\kappa$-concentrated} on a countable set $D\sub X$ if
$\card{X\sm U}<\kappa$ for all open sets $U$ containing $D$.

Every compact set $K\sub\roth$ is $\le$-bounded. A classic argument of 
Lawrence~\cite[Propositions~2--3]{Law} 
implies that, for each cardinal number $\kappa$, the existence of a $\kappa$-concentrated real 
set is equivalent to the existence of a $\kappa$-unbounded set in $\roth$. Essentially, this is due
to the following fact.

\blem\label{lem:UnbImplCon}
Let $\kappa$ be a cardinal number, and $X\sub\roth$ be a set with $\card{X}\ge\kappa$.
The set $X$ is $\kappa$-unbounded if and only if the real set $X\cup\Fin$
is $\kappa$-concentrated on $\Fin$.
\elem
\bpf
$(\Rightarrow)$
Let $U\sub\PN$ be an open set containing the set $\Fin$.
The set $K:=\PN\sm U$ is a closed, and thus compact, subset of $\PN$.
Since $U\spst\Fin$, we have $K\sub\roth$.
Since compact subsets of $\roth$ are $\le$-bounded and
the set $X$ is $\kappa$-unbounded, we have 
\[
\card{(X\cup\Fin)\cap K}=\card{X\cap K}<\kappa.
\]

$(\Leftarrow)$ For each bound $b\in\roth$, the set $K:=\sset{a\in\roth}{a\le b}$ is compact.
Thus, the set $U:=\PN\sm K$ is an open set containing $\Fin$, and
we have $\card{X\sm U}<\kappa$.
\epf

A set $X\sub\roth$ is \emph{dominating} if for each function $a\in\roth$ there is
a function $x\in X$ such that $a\les x$. Let $\fd$ be the minimal cardinality of a dominating
set in $\roth$. Much information about the cardinal number $\fd$, and about other ones defined below,
is available~\cite{BlassHBK}.
Every real set of cardinality smaller than $\fd$ is Menger, and no dominating
subset of $\roth$ is Menger~\cite[Theorem~4.4]{coc2}.
The former assertion implies that every $\fd$-concentrated real set is
Menger.\footnote{Moreover, $\fd$-concentrated sets have the stronger selective property
$\sone(\Ga,\Op)$~\cite{Ideals,MHP}.}

\bcor\label{cor:dUnbIsM}
For each $\fd$-unbounded set $X\sub\roth$,
the real set $X\cup\Fin$ is Menger.\qed
\ecor

\bdfn
For functions $a,b\in\roth$, we write $a\lei b$ if $b\nless^* a$, that is,
if $a(n)\le b(n)$ for infinitely many natural numbers $n$.
For a set $X\sub\roth$ and a function $b\in\roth$,
we write $X\lei b$ if $x\lei b$ for each function $x\in X$.
This convention applies to all binary relations.
\edfn

There are, provably, $\fd$-unbounded sets and $\cf(\fd)$-unbounded sets:
Let $\sset{d_\alpha}{\alpha<\fd}$ be a dominating set.
For each ordinal number $\alpha<\fd$,
take a function $x_\alpha\in\roth$ such that $\sset{d_\beta,x_\beta}{\beta<\alpha}<^\oo x_\alpha$.
Then the set $\sset{x_\alpha}{\alpha<\fd}$ is $\fd$-unbounded.
Taking a cofinal subset $I\sub\fd$ of cardinality $\cf(\fd)$, we obtain
the $\cf(\fd)$-unbounded set $\sset{x_\alpha}{\alpha\in I}$.

\blem\label{lem:ergosum}
For sets $a,b\in\PN$, let
\[
a\uplus b:=(2a)\cup(2b+1)=\set{2k}{k\in a}\cup\set{2k+1}{k\in b}.
\]
Then:
\be
\item For each set $a\in\roth$ and each natural number $n$, we have 
$(a\uplus a)(2n) =2a(n)+1$.
\item For all sets $a,b,c,d\in\roth$ with $a\le b$ and $c\le d$, we have $a\uplus c\le b\uplus d$.\qed
\ee
\elem

\bthm\label{thm:master}
Let $\kappa\in\{\cf(\fd),\fd\}$, and
$X\sub\roth$ be a set containing a $\kappa$-unbounded set.
There is a
$\fd$-concentrated real set $Y$
such that the planar set $X\x Y$ is not Menger.
\ethm
\bpf
Let $A\sub X$ be a $\kappa$-unbounded set. By moving to
a subset of $A$, we may assume that $\card{A}=\kappa$.
Let $D\sub\roth$ be a dominating set of cardinality $\fd$.
Decompose
\[
D=\Un_{a\in A}I_a
\]
such that $\card{\Un_{a\in B}I_a}<\fd$ for all sets $B\sub A$
of cardinality smaller than $\kappa$.
(If $\kappa=\fd$, we can take every set $I_a$ to be a singleton).
Fix elements $a\in A$ and $d\in I_a$. Take a function $d'\in\roth$ such that
$a,d\le d'$.
Consider the set $\sset{a\uplus d'}{a\in A, d\in I_a}$.
Its cardinality is at most $\fd$, and since its projection on the odd coordinates is dominating,
its cardinality is exactly $\fd$.

\bclm
The set $\sset{a\uplus d'}{a\in A, d\in I_a}$ is $\fd$-unbounded.
\eclm
\bpf
Let $b\in\roth$. Define $b'(n):=b(2n)$ for all natural numbers $n$.
Let $K:=\sset{a\in A}{a\le b'}$. Then $\card{K}<\kappa$.

Let $a\in A\sm K$ and $d\in I_a$. There is a natural number $n$ such that
\[
b(2n)=b'(n) < a(n)\le  2a(n)+1=(a\uplus a)(2n)\le (a\uplus d')(2n),
\]
and thus $a\uplus d'\nleq b$.
Therefore,
\[
\card{\set{a\uplus d'}{a\in A, d\in I_a, a\uplus d'\le b}}\le\card{\set{a\uplus d'}{a\in K, d\in I_a}}<\fd.\qedhere
\]
\epf

By Lemma~\ref{lem:UnbImplCon},
the real set
\[
Y:=\set{a\uplus d'}{a\in A, d\in I_a}\cup\Fin
\]
is $\fd$-concentrated on the set $\Fin$. In particular, the set $Y$ is Menger.

For sets $a,b\in\PN$, let $a\oplus b$ denote the symmetric difference of the sets $a$ and $b$.
With respect to the operator $\oplus$, the space $\PN$ is a topological group.

\bclm
The set $(2X)\oplus Y$ is a dominating subset of $\roth$.
\eclm
\bpf
For all sets $a,b,c\in\PN$, we have
$(2a)\oplus (b\uplus c)=(a\oplus b)\uplus c\spst 2c+1$.
It follows that $(2X)\oplus Y\sub\roth$.

For each function $d$ in the dominating set $D$ we started with,
let $a\in A$ be a function such that $d\in I_a$.
As $a\in X$ and $a\uplus d'\in Y$, we have 
\[
2a\oplus (a\uplus d')=(a\oplus a)\uplus d'=\emptyset\uplus d'=2d'+1\in (2X)\oplus Y.
\]
Since $d\le d'\le 2d'+1$ for all functions $d\in D$, the set $(2X)\oplus Y$ is dominating.
\epf

In summary, the set $(2X)\oplus Y$ is a continuous image of the planar set $X\x Y$ in $\roth$
that is dominating. It follows that the space  $X\x Y$ is not Menger.
\epf

Let $X$ be a real set of cardinality smaller than $\fd$.
Then the set $X$ is \emph{trivially} Menger: the topology used is irrelevant,
as long as we restrict attention to countable covers.
In particular, all finite powers of the real set $X$ are Menger, even for
countable \emph{Borel} covers; a strong property~\cite{cbc}.

\bthm\label{thm:dSing}
Assume that $\cf(\fd)<\fd$. There are real sets $X$ and $Y$ such that
$\card{X}<\fd$ and the set $Y$ is $\fd$-concentrated, but the planar set
$X\x Y$ is not Menger.
\ethm
\bpf
By the discussion preceding Lemma~\ref{lem:ergosum}, there are
$\cf(\fd)$-unbounded sets in $\roth$.
Apply Theorem~\ref{thm:master} to any of these sets.
\epf

Let $\kappa$ be a cardinal number. A real set of cardinality at least $\kappa$ is 
\emph{$\kappa$-Luzin} if the cardinalities of its intersections with 
meager sets are all smaller than $\kappa$.
Let $\cov(\cM)$ be the minimal cardinality of a cover of the real line
by meager sets,
and $\cof(\cM)$ be the minimal cardinality of a cofinal family of meager real sets.
The hypothesis $\cov(\cM)=\cof(\cM)$ implies that
there are $\cov(\cM)$-Luzin sets whose product is not Menger~\cite[Theorem~49]{cbc}.
Since $\cov(\cM)\le\fd$, every $\cov(\cM)$-Luzin set is $\fd$-concentrated, and thus Menger.
In general, $\cov(\cM)\le\fd\le\cof(\cM)$, and thus
the following corollary implies (using the same hypothesis) that
for \emph{every} $\cov(\cM)$-Luzin set $L$ there is a $\fd$-concentrated real set
$Y$ such that the planar set $L\x Y$ is not Menger.

\bcor\label{cor:LnonpM}
Let $\kappa\in\{\cf(\fd),\fd\}$. For each
$\kappa$-Luzin set $L$, there is a $\fd$-concentrated real set $Y$ such that
the planar set $L\x Y$ is not Menger.
In particular, if $\aleph_1=\cf(\fd)$, then this is the case for every Luzin set.
\ecor
\bpf
By applying a homeomorphism, we may assume that $L\sub\roth$.
Every $\kappa$-Luzin subset of $\roth$ is $\kappa$-unbounded.
Apply Theorem~\ref{thm:master}.
\epf

The most important application of Theorem~\ref{thm:master} appears in the next section.

\section{Bi-$\fd$-unbounded sets}

For a set $a\in\PN$, let $a\comp:=\bbN\sm a$.
Let $\ici:=\sset{a\in\roth}{a\comp\in\roth}$, the family of infinite co-infinite subsets
of $\bbN$.

\bdfn
Let $\kappa$ be an infinite cardinal number.
A set $X\sub\ici$ is \emph{bi-$\kappa$-unbounded}
if the sets $X$ and $\sset{x\comp}{x\in X}\sub\roth$ are both $\kappa$-unbounded.
\edfn

\bthm
\label{thm:bidiNotpM}
Let $\kappa\in\{\cf(\fd),\fd\}$.
Let $X\sub\roth$ be a bi-$\kappa$-unbounded set. Then:
\be
\item The real set $X\cup\Fin$ is $\kappa$-concentrated.
In particular, it is Menger.
\item There is a $\fd$-concentrated real set $Y$ such that the
planar set $(X\cup\Fin)\x Y$ is not Menger.
\ee
\ethm
\bpf
(1) By Corollary~\ref{cor:dUnbIsM}.

(2) The continuous image $\sset{x\comp}{x\in X\cup\Fin}$ of the real set $X\cup\Fin$ in $\PN$ is
a $\kappa$-unbounded subset of $\roth$. Apply Theorem~\ref{thm:master}.
\epf

The existence of bi-$\fd$-unbounded sets and bi-$\cf(\fd)$-unbounded sets
is a mild hypothesis.
A set $r\in\roth$ \emph{reaps} a family $A\sub\roth$ if, for each set $a\in A$,
both sets $a\cap r$ and $a\sm r$ are infinite.
Let $\fr$ be the minimal cardinality of a family $A\sub\roth$ that no
set $r$ reaps.

For natural numbers $n<m$, let $[n,m):=\{n,n+1,\dots,m-1\}$.

\bthm\label{thm:ExistsbiDunb}
The following assertions are equivalent:
\be
\item $\fd\le\fr$.
\item There are bi-$\fd$-unbounded sets in $\roth$.
\item There are bi-$\cf(\fd)$-unbounded sets in $\roth$.
\ee
\ethm
\bpf
$(1)\Impl (2),(3)$: We use the following lemma, to which we provide a short, direct proof.

\blem[Mej\'\i{}a~\cite{Mejia}]
\label{lem:mej}
Let $X\sub\roth$. If $\card{X}<\min\{\fd,\fr\}$, then there is an element
$b\in\ici$ such that $X\leinf b$ and $X\leinf b\comp$.
\elem
\bpf
For a set $x\in\roth$ with $1\notin x$,
define a function $\tilde x\in\roth$ by
$\tilde x(1):=x(1)$, and $\tilde x(n+1):= x(\tilde x(n))$ for each natural number $n$.

We may assume that $1\notin x$ for all sets $x\in X$.
Since $\card{X}<\fd$, there is
a function $a\in\roth$ such that the sets
\[
I_x := \set{n}{\card{[a(n),a(n+1))\cap \tilde x}\ge 2}
\]
are infinite for all sets $x\in X$~\cite[Theorem~2.10]{BlassHBK}.
Since $\card{X}<\fr$, there is a set $r\in\roth$ that reaps
the family $\set{I_x}{x\in X}$. Define
\[
b:=\Un_{n\in r}[a(n),a(n+1)).
\]
Fix a set $x\in X$. Let $n$ be a member of the infinite set $r\cap I_x$.
There are at least two elements in the set
$[a(n),a(n+1))\cap \tilde x$; let $\tilde x(i)$ be the minimal one.
Then $x(\tilde x(i))=\tilde x(i+1)\in [a(n),a(n+1))$.
Since $n\in r$, the set $b\comp\cap [a(n),a(n+1))$ is empty, and thus
$a(n+1)\le b\comp(\tilde x(i))$. It follows that
$x(\tilde x(i))<b\comp(\tilde x(i))$.
Similarly, every number $n\in I_x\sm r$ produces a number $i$
such that $x(\tilde x(i))<b(\tilde x(i))$.
\epf

Let $\sset{d_\alpha}{\alpha<\fd}\sub\roth$ be a dominating set.
By Lemma~\ref{lem:mej}, for each ordinal number $\alpha<\fd$,
there is a set $x_\alpha\in\ici$ such that
$\sset{d_\beta,x_\beta}{\beta<\alpha}<^\oo x_\alpha,x_\alpha\comp$.
Then the set $\sset{x_\alpha}{\alpha<\fd}$ is bi-$\fd$-unbounded.
Let $I$ be a cofinal subset of the cardinal number $\fd$, of cardinality $\cf(\fd)$.
Then the set $\sset{x_\alpha}{\alpha\in I}$ is bi-$\cf(\fd)$-unbounded.

$(2)\Impl(1)$:
We may assume that the cardinal number $\fd$ is regular.
Indeed, it is known that if $\fr<\fd$ then $\fd$ is regular.\footnote{In the notation of Section~\ref{sec:FM}, fix an
ultrafilter $U$ with pseudobase $P$ of cardinality $\fr$~\cite[Theorem~9.9]{BlassHBK}, and
take a $\le_U$-dominating set $D$ of cardinality $\bof(U)$, a regular cardinal number.
Then the set $\set{f\circ p}{f\in D, p\in P}$
is dominating, and thus $\fd\le\bof(U)(\le\fd)$.}
Thus, if $\fd$ is singular, then $\fd\le\fr$, and we are done.

Let $X\sub\roth$ be a bi-$\fd$-unbounded set.
Let $A\sub\roth$ be a family with $\card{A}<\fd$.
We prove that the family $A$ is reapable.
We may assume that for each set $a\in A$ and each finite set $s$,
we have $a\sm s\in A$.

Since the set $X$ is bi-$\fd$-unbounded, the set
\[
\Un_{a\in A}\set{x\in X}{x\les a\mbox{ or }x\les a\comp}
\]
is a union of less than $\fd$ sets, each of cardinality smaller than $\fd$.
Thus, there is an element $r\in X$ that is not included in that set,
that is, such that $A<^\oo r,r\comp$.

The set $r$ reaps the family $A$:
Fix a set $a\in A$.
Assume that the set $a\cap r$ is finite.
Then the set $a':=a\sm r$ is in $A$, and thus $a'<^\oo r\comp$.
But $a'\sub r\comp$, and thus $r\comp\le a'$; a contradiction.
For the same reason, the set $a\sm r$ is infinite, too.

$(3)\Impl(1)$: If the cardinal number $\fd$ is regular, then the
previously established implication applies. And if it is singular, then as explained in
the previous implication, we have $\fd\le\fr$. In either case, we are done.
\epf

A topological space $X$ is \emph{Rothberger} if for each sequence
$\eseq{\cU}$
of open covers of $X$, there are elements
$U_1\in\cU_1$, $U_2\in\cU_2$, \dots
with $X\sub\Un_nU_n$.
Every real set of cardinality smaller than $\cov(\cM)$ is Rothberger~\cite[Theorem~4.2]{coc2},
and therefore so is every $\cov(\cM)$-concentrated real set. Since $\cov(\cM)\le\fr$~\cite[Theorem~5.19]{BlassHBK}, 
we obtain the following result.

\bcor\label{cor:Roth}
Assume that $\cov(\cM)=\fd$.
Then there are two Rothberger real sets whose product is not Menger.
\qed
\ecor

\section{Filter-Menger spaces}
\label{sec:FM}

For sets $a,b\in\roth$, we write $a\as b$ if the set $a\sm b$ is finite.
A \emph{semifilter}~\cite{BZCSFQdM} is a set $S\sub\roth$ such that, for each set
$s\in S$ and each set $b\in\roth$ with $s\as b$, we have $b\in S$.\footnote{Semifilters
are normally denoted by calligraphic letters.
Here, we view them as sets of points in, and thus subspaces of, the Cantor space $\PN$.
Thus, we use the standard typefaces, as we do for arbitrary points and sets in topological spaces.}
Important examples of semifilters include the maximal semifilter $\roth$, the minimal semifilter
$\cFin$ of all cofinite sets, and every nonprincipal ultrafilter on $\bbN$.

Let $S$ be a semifilter. For functions $a,b\in\roth$, let
\[
[a\le b] := \set{n}{a(n)\le b(n)}.
\]
We write $a\leS b$ if $[a\le b]\in S$.
Let $\bof(S)$ be the minimal cardinality of a $\leS$-unbounded subset of $\roth$.
For a semifilter $S$, let $S^+:=\sset{a\in\roth}{a\comp\notin S}$.
For all sets $a\in S$ and $b\in S^+$, the intersection $a\cap b$ is infinite.
For functions $a,b\in\roth$, we have that $a\nleq_S b$ if and only if $b<_{S^+} a$.
The $\kappa$-unbounded sets presented in the previous sections
are instances of the following notion, which generalizes
the earlier notion of $\bof(S)$-scale~\cite[Definition~2.8]{sfh}.

\bdfn
Let $S$ be a semifilter. A set $X\sub\roth$ with $\card{X}\ge\bof(S)$
is an \emph{$S$-scale} if, for each function $b\in\roth$, there is
a function $c\in\roth$ such that
\[
b\le_{S^+} c\le_S x
\]
for all but less than $\bof(S)$ functions $x\in X$.
\edfn

\bprp[{\cite[Lemma~2.9]{sfh}}]\label{prp:ExsistScale}
For each semifilter $S$ there is an $S$-scale.
\eprp
\bpf
Let $\sset{b_\alpha}{\alpha<\bof(S)}\sub\roth$ be a $\leS$-unbounded set.
For each ordinal number $\alpha<\bof(S)$,
there is a function $x_\alpha\in\roth$ such that
$\sset{b_\beta,x_\beta}{\beta<\alpha}<_S x_\alpha$. 
The set $\sset{x_\alpha}{\alpha<\bof(S)}$ is an $S$-scale.
Indeed, fix a function $b\in\roth$.
There is an ordinal number $\beta<\bof(S)$ such that $b_\beta\nleq_S b$,
and thus $b\le_{S^+} b_\beta$. For each ordinal number $\alpha>\beta$, we have
$b_\beta\leS x_\alpha$.
\qedhere
\epf

Let $S$ be a semifilter, and $b,c,x\in\roth$ functions satisfying $b\le_{S^+} c\le_S x$.
Then the set $[b\le x]$ contains the intersection $[b\le c]\cap [c\le x]$ of an element
of $S^+$ and an element of $S$. In particular, we have $b\leinf x$.

\bprp\label{prp:bofS-con}
Let $S$ be a semifilter.
Every $S$-scale is a $\bof(S)$-unbounded subset of $\roth$, and thus
its union with the
set $\Fin$ is $\bof(S)$-concentrated. In particular,
no union of an $S$-scale and $\Fin$ is $\sigma$-compact~\cite[Lemma~1.6]{MHP}.
\qed
\eprp

Let $S$ be a semifilter. A topological space $X$ is \emph{$S$-Menger} if
for each sequence
$\eseq{\cU}$
of open covers of $X$, there are finite sets
$\cF_1\sub\cU_1$, $\cF_2\sub\cU_2$, \dots
such that $\sset{n}{x\in\Un\cF_n}\in S$ for all points $x\in X$.
A topological space is Menger if and only if it is $\roth$-Menger.
For the filter $\cFin$ of cofinite sets,
the property $\cFin$-Menger is the classic \emph{Hurewicz property}~\cite{Hure25}.
Thus, for every semifilter $S$, we have the following implications.
\[
\mbox{Hurewicz} \Longrightarrow S\mbox{-Menger} \Longrightarrow \mbox{Menger}.
\]

A function $\Psi$ from a topological space $X$ into $\roth$ is \emph{upper continuous} if the sets
$\sset{x\in X}{\Psi(x)(n)\le m}$ are open for all natural numbers $n$ and $m$.
In particular, continuous functions are upper continuous.
By earlier methods~\cite[Theorem~7.3]{SPMProd}, we have the following result.

\bprp\label{prp:SMenChar}
Let $X$ be a topological space, and $S$ be a semifilter.
The following assertions are equivalent:
\be
\itm The space $X$ is $S$-Menger.
\itm The space $X$ is Lindel\"of, and every upper continuous image of $X$ in $\roth$
is $\leS$-bounded.\qed
\ee
\eprp

For especially nice classes of spaces,
such as Lindel\"of zero-dimensional spaces or real sets,
\emph{upper continuous} can be replaced by \emph{continuous}
in Proposition~\ref{prp:SMenChar}.
In general, however, this is not the case:
The properties considered here are hereditary for closed subsets.
Consider the planar set
\[
X := ((\bbR\sm\bbQ)\x [0,1])\cup (\bbR\x \{1\})\sub\bbR^2.
\]
This set is not Menger,
since the non-Menger set $(\bbR\sm\bbQ)\x\{0\}$ (homeomorphic to $\roth$) is closed in $X$.
Since the set $X$ is connected, every \emph{continuous} image of $X$ in
$\roth$ is a singleton.

For a set $a\in S^+$, let
\[
S_a:=\set{c\in\roth}{\exists s\in S, s\cap a\as c},
\]
the semifilter generated by the sets $\sset{s\cap a}{s\in S}$.
The following observation generalizes an earlier result~\cite[Theorem~2.14]{sfh}.

\bprp\label{prp:SscaleSetBdd}
Let $S$ be a semifilter, and $X\sub\roth$ be an $S$-scale.
Every upper continuous image of the real set $X\cup\Fin$ in $\roth$ is
$\le_{S_a}$-bounded for some set $a\in S^+$.
\eprp
\bpf
Let $\Psi\colon X\cup\Fin\to\roth$ be an upper continuous function.
We use the forthcoming Lemma~\ref{lem:PTWL},
in the case $Y=\{0\}$. This special case
was, implicitly, established by Bartoszy\'nski and Shelah~\cite[Lemma~2]{BaSh01}.
This lemma provides a function $b\in\roth$ such that $\Psi(x)(n)\le b(n)$ for all
functions $x\in X$ and all natural numbers $n$ with
$b(n)\le x(n)$.
Since the set $X$ is an $S$-scale, there is a function $c\in\roth$ such that
$b\le_{S^+} c\leS x$ for all but less than $\bof(S)$ functions $x\in X$.
For these points $x$, we have that $\Psi(x)(n)\le b(n)$ for all natural numbers
$n\in [b\le c]\cap [c\le x]$. Take $a:=[b\le c]$.

The image of the remaining points of the set $X\cup\Fin$ is $\leS$-bounded
by some member $b'\in\roth$. Then any function $d\in\roth$ with $b,b'\le d$
is a bound as required.
\epf

By \emph{filter} we mean a semifilter closed under finite intersections.
If $F$ is a filter, then $a\cap b\in F^+$ for all sets $a\in F$
and $b\in F^+$. And if $F$ is an ultrafilter, then $F^+=F$.

\bcor\label{cor:simple}
For every filter $F$, the union of every $F$-scale and $\Fin$ is $F^+$-Menger,
and if $F$ is an ultrafilter, this union is $F$-Menger.\qed
\ecor

Let $\fb$ be the minimal cardinality of a $\les$-unbounded subset of $\roth$.

\bthm\label{thm:strongbidiscale}
Assume that $\fb=\fd$. Let $S$ be a semifilter. The following assertions are equivalent:
\be
\item The semifilter $S$ is nonmeager.
\item There are an $S$-scale $X\sub\roth$ and a
$\fd$-concentrated real set $Y$ such that the planar set $(X\cup\Fin)\x Y$ is not Menger.
\ee
\ethm
\bpf
$(1)\Impl (2)$:
Let $\sset{d_\alpha}{\alpha<\fd}$ be a dominating subset of $\roth$. Fix an ordinal number $\alpha<\fd$.
Since $\fb=\fd$, there is a function $b\in\roth$ such that
$\sset{d_\beta,x_\beta}{\beta<\alpha}<^* b$.
The set $\sset{x\in\ici}{b\lei x\comp}$ is comeager.
Since the semifilter $S$ is nonmeager, the set
$\sset{x\in\ici}{b\leS x}$ is nonmeager~\cite[Corollary~3.4]{sfh}.
Thus, there is a set $x_\alpha\in\ici$ such that
$b\leS x_\alpha$ and $b\lei x_\alpha\comp$.
Then the set $X:=\sset{x_\alpha}{\alpha<\fd}$ is an $S$-scale, and
it is bi-$\fd$-unbounded. Apply Theorem~\ref{thm:master}.

$(2)\Impl (1)$:
Let $S$ be a meager semifilter, and $X\sub\roth$ be an $S$-scale. 
By Theorem~\ref{thm:HH} below, the real set $X\cup\Fin$ is, in particular, Hurewicz.
Products of Hurewicz sets and $\fd$-concentrated real sets are Menger~\cite[Theorem~4.6]{AddGN}.
\epf

The real set $X\cup\Fin$ in Theorem~\ref{thm:strongbidiscale}
is not Hurewicz since its image under the function $x\mapsto x\comp$ is unbounded in $\roth$.
The existence of non-Hurewicz sets of this form follows
from a weaker hypothesis~\cite[Theorem 3.9]{sfh},
but without the non-productive property.
The product of every Hurewicz real set and every $\fd$-concentrated real set
is Menger~\cite[Theorem 4.6]{AddGN}. The following theorem implies that
this assertion cannot be established for spaces that are not Hurewicz.

Let P be a property of topological spaces. A real set $X$ is \emph{productively P} if
for each topological space $Y$ with the property P, the product space $X\x Y$ has the property P.
The question whether productively P implies productively Q, for P and Q covering properties among
those studied here, has a long history. The remainder of this paragraph assumes 
that $\fd=\aleph_1$.
Aurichi and Tall~\cite{AurTall} improved several earlier results by proving that
every productively Lindel\"of space is Hurewicz. It was later shown that every 
productively Lindel\"of space is productively Hurewicz and productively Menger~\cite[Theorem~8.2]{SPMProd}. Thus, productively Lindel\"of implies productively Menger, and the following theorem shows that
productively Menger suffices to imply productively Hurewicz.

\bthm\label{thm:pMispH}
Assume that $\fb=\fd$.
\be
\item For every unbounded set $X\sub\roth$, there is a $\fd$-concentrated real set
$Y$ such that the planar set $X\x Y$ is not Menger.
\item In the realm of hereditarily Lindel\"of spaces: If a real set $X$ is
productively Menger, then it is productively Hurewicz.
\ee
\ethm
\bpf
(1) Let $\set{d_\alpha}{\alpha<\fd}$ be a dominating set in $\roth$.
Since $\fb=\fd$,
for each ordinal number $\alpha<\fd$ the set $\set{d_\beta}{\beta<\alpha}$ is
bounded, and thus there is a function $x_\alpha\in X$ such that
$\set{d_\beta,x_\beta}{\beta<\alpha}<^\oo x_\alpha$. 
Then the subset $\set{x_\alpha}{\alpha<\fd}$ of the set
$X$ is $\fd$-unbounded, and Theorem~\ref{thm:master} applies.

(2) Assume that there is a Hurewicz hereditarily Lindel\"of space $H$ such that
the product space $X\x H$ is not Hurewicz.
Then there is an unbounded upper continuous image $Z$ of the space
$X\x H$ in $\roth$~\cite[Theorem~7.3]{SPMProd}.
By (1), there is a $\fd$-concentrated real set $Y$ such that
the planar set $Z\x Y$ is not Menger. Since the set $Z\x Y$ is an upper continuous
image of the product space $X\x H\x Y$, the latter space is not Menger, too.
As the space $H$ is Hurewicz and hereditarily Lindel\"of and the set $Y$ is
$\fd$-concentrated, the product space $H\x Y$ is Menger~\cite[Theorem 4.6]{AddGN}.
In summary, the product of the real set $X$ and the Menger, hereditarily Lindel\"of space $H\x Y$
is not Menger.
\epf

Some special hypothesis is necessary for Theorem~\ref{thm:pMispH}:
The union of less than $\fg$ Menger real sets is Menger~\cite{SF1, sft}.
Assume that $\fb<\fg$. Then any unbounded real set $X\sub\roth$ of cardinality $\fb$
is productively Menger but not Hurewicz.

\section{Productive real sets}
\label{sec:pR}

In this section, we establish preservation of some properties under products.
We begin with a generalization of an earlier result~\cite[Lemma 6.3]{SPMProd}
to general topological spaces. The earlier proof~\cite[Lemma 6.3]{SPMProd}
does not apply in this general setting; we provide an alternative proof.

\blem[Productive Two Worlds Lemma]\label{lem:PTWL}
Let $X$ be a subset of $\roth$, $Y$ be an arbitrary space,
and $\Psi\colon (X\cup\Fin)\times Y\to\roth$ be an
upper continuous function.
There is an upper continuous function $\Phi\colon Y\to\roth$ such that, for all
points $x\in X$ and $y\in Y$, and all natural numbers $n$:
\begin{center}
If  $\Phi(y)(n)\le x(n)$, then $\Psi(x,y)(n)\le\Phi(y)(n)$.
\end{center}	
\elem
\bpf
For natural numbers $n$ and $m$, let $U^n_m:=\Psi\inv[\sset{a\in\roth}{a(n)\le m}]$.
For each natural number $n$, the family $\sset{U^n_m}{m\in\bbN}$
is an ascending open cover of the product space $(X\cup\Fin)\x  Y$.
By enlarging the sets $U^n_m$, we may assume that they are open in the
larger space $\PN\times Y$.

Fix a point $y\in Y$ and a natural number $n$.
Set $a_y(1):=1$ and $V^y_1:=Y$.
For a natural number $k$,
let $m$ be the minimal natural number with
$\PS([1,a_y(k)))\x  \{y\}\sub {U}^n_m$.
Let $a_y(k+1)$ be the minimal natural number such that 
$a_y(k+1)\ge m$ and $\PS([a_y(k),a_y(k+1))\comp)\x  \{y\}\sub {U}^n_m$.
Since our open covers are ascending, we have
\[
\PS([a_y(k),a_y(k+1))\comp)\x  \{y\}\sub {U}^n_{a_y(k+1)}.
\]
Notice that the number $a_y(k+1)$ is minimal with this property.
As the set $\PS([a_y(k),a_y(k+1))\comp)\x  \{y\}$ is compact,
there is an open neighborhood
$V^y_{k+1}\sub V^y_{k}$ of the point $y$ such that
\[
\PS([a_y(k),a_y(k+1))\comp)\x  V^y_{k+1}\sub {U}^n_{a_y(k+1)}.
\]
Define
\[
\Phi(y)(n):=a_y(n+1).
\]
For each point $y'\in V^y_{n+1}$, we have $y'\in V^y_k$ for all $k=1,\dots,n+1$.
The sequence $a_{y'}(1),\dots,a_{y'}(n+1)$
is bounded by the sequence $a_y(1),\dots,a_y(n+1)$, coordinate-wise:
For $k=1$, we have $a_{y'}(1)=a_{y}(1)=1$. Assume that 
$a_{y'}(k)\le a_{y}(k)$. Then
\[
\PS([a_{y'}(k),a_y(k+1))\comp)\x  \{y'\}\sub \PS([a_y(k),a_y(k+1))\comp)\x  V^y_{k+1}\sub {U}^n_{a_y(k+1)},
\]
and, by the minimality of the number $a_{y'}(k+1)$,
we have $a_{y'}(k+1)\le a_{y}(k+1)$, too.

In summary, we have $\Phi(y')(n)=a_{y'}(n+1)\le a_y(n+1)\le\Phi(y)(n)$
for all points $y'\in V^y_{n+1}$.
This shows that the function $\Phi$ is upper continuous.

Fix a point $y\in Y$ and a natural number $n$.
Let $x\in\roth$ be an element with $\Phi(y)(n)\le x(n)$.
As $a_y(n+1)=\Phi(y)(n)\le x(n)$, there is a natural number $k\le n$ such that
$x\cap [a_y(k),a_y(k+1))=\emptyset$.
Thus,
\[
(x,y)\in \PS([a_y(k),a_y(k+1))\comp)\x  \{y\} \sub U^n_{a_y(k+1)}\sub U^n_{\Phi(y)(n)},
\]
and therefore $\Psi(x,y)(n)\le\Phi(y)(n)$.
\epf

For filters, we obtain a productive version of Proposition~\ref{prp:SscaleSetBdd}.

\bthm\label{thm:master2}
Let $F$ be a filter and $X\sub\roth$ be an $F$-scale.
For each $F$-Menger space $Y$,
every upper continuous image of the product space $(X\cup\Fin)\x Y$ in $\roth$ is
$\le_{F_a}$-bounded for some set $a\in F^+$.
\ethm
\bpf
Let $\Psi\colon (X\cup\Fin)\x Y\to\roth$ be an upper continuous function.
Let $\Phi\colon Y\to\roth$ be as in the
Productive Two Worlds Lemma (Lemma~\ref{lem:PTWL}).
Since the space $Y$ is $F$-Menger, there is a function $b\in\roth$
such that $\Phi[Y]\le_F b$.
As the set $X$ is an $F$-scale, there is a function $c\in\roth$
such that $b\le_{F^+} c\le_F x$ for all but less than
$\bof(F)$ elements of $X$. Let $a:=[b\le c]$, an element of the semifilter $F^+$.
Then the cardinality of the set
\[
Z:=\set{x\in X}{b\nleq_{F_{a}} x}
\]
is smaller than $\bof(F)$.

Fix a pair $(x,y)\in (X\sm Z)\x Y$. Then $b\le_{F_{a}} x$ and $\Phi(y)\le_F b$.
Since $F$ is a filter, we have 
\[
[\Phi(y)\le x]\spst [\Phi(y)\le b]\cap [b\le x]\in F_a.
\]
This shows that $\Psi[(X\sm Z)\x Y]\le_{F_a} b$.
Let $z\in Z\cup\Fin$. Since
the set $\{z\}\x Y$ is $F$-Menger,
$\Psi[\{z\}\x Y]\le_{F} c_z$  for some function $c_z\in\roth$. Since $\card{Z\cup\Fin}<\bof(F)$,
there is a function $c\in\roth$ such that $\sset{c_z}{z\in Z\cup\Fin}\le_F c$.
As $F$ is a filter, we have 
$\Psi[(Z\cup\Fin)\x Y]\le_{F} c$,
and therefore
$\Psi[(X\cup\Fin)\x Y]\le_{F_a} \set{\max\{b(n),c(n)\}}{n\in\bbN}$.
\epf

\bthm\label{thm:UU}
Let $F$ be a filter, and $X\sub\roth$ be an $F$-scale. In the realm of hereditarily Lindel\"of spaces:
\be
\item For each $F$-Menger space $Y$, the product space $(X\cup\Fin)\x Y$ is $F^+$-Menger.
\item If $F$ is an ultrafilter, then the real set $X\cup\Fin$ is productively $F$-Menger.
\ee
\ethm
\bpf
Every product of a metrizable Lindel\"of space and a
hereditarily Lindel\"of spaces is Lindel\"of.
Apply Theorem~\ref{thm:master2}.
\epf

The following theorem was previously known for $\fb$-scales,
a special kind of $\cFin$-scales~\cite[Theorem~6.5]{SPMProd}. 
This theorem and the subsequent one
improve upon  earlier results~\cite[Corollary~4.4]{sfh}, asserting that the
corresponding properties hold in all finite powers.

A semifilter $S$ is meager if and only
if there is a function $h\in\roth$ such that for each set $s\in S$, the set $s\cap[h(n),h(n+1))$ is
nonempty for almost all natural numbers $n$~\cite[Theorem~21]{Talagrand}.
For meager semifilters $S$, we have $\bof(S)=\fb$~\cite[Corollary~2.27]{sfh},
and $S$-Menger is equivalent to Hurewicz~\cite[Theorem~2.32]{sfh}. The following theorem
generalizes an earlier result~\cite[Theorem~2.28]{sfh}, using a similar proof.

\bthm\label{thm:HH}
Let $S$ be a meager semifilter, and $X\sub\roth$ be an $S$-scale. Then, in the realm of
hereditarily Lindel\"of spaces, the real set $X\cup\Fin$ is productively Hurewicz.
\ethm
\bpf
Let $Y$ be a hereditarily Lindel\"of, Hurewicz space.
Since the space $Y$ is hereditarily Lindel\"of, the product space $(X\cup\Fin)\x Y$ is Lindel\"of.
Let $\Phi\colon X\x Y\to\roth$ be an upper continuous function and
$\Psi\colon Y\to \roth$ be the upper continuous function provided by
the Productive Two Worlds Lemma (Lemma~\ref{lem:PTWL}).
Since the space $Y$ is Hurewicz, its image $\Psi[Y]$ is $\les$-bounded
by some function $b\in\roth$.
Let $h\in\roth$ be a witness for the semifilter $S$ being meager.
Define a function $\tilde{b}\in\NN$ by
\[
\tilde{b}(k):=b(h(n+2))
\]
for all natural numbers $n$ and for $k\in [h(n),h(n+1))$.
Then  $\Psi[Y]\les b\le \tilde{b}$.

Since the set $X$ is an $S$-scale, there is a function $c\in\roth$ such that 	$\tilde{b}\le_{S^+} c$ and all but less than $\fb$ functions $x\in X$ belong to the set
\[
\tilde{X}:=\set{x\in X}{c\le_S x}.
\]

\bclm
The set $\Phi[(\tilde{X}\cup\Fin)\x Y]$ is $\les$-bounded.	
\eclm
\bpf
Fix  a function $x\in\tilde{X}$. Then $[c\le x]\in S$, and thus
the set $[c\le x]\cap [h(n),h(n+1))$ is nonempty for
almost all natural numbers $n$.
Let
\[
d:=\set{n\in\bbN}{[\tilde{b}\le c]\cap[h(n-1),h(n))\neq\emptyset}.
\]
Then, for almost all natural numbers $n\in d$,
there are natural numbers $l\in [\tilde{b}\le c]\cap[h(n-1),h(n))$ and
$m\in[c\le x]\cap [h(n),h(n+1))$, and we have
\[
b(h(n+1))=\tilde{b}(l)\le c(l)\le c(m)\le x(m)\le x(h(n+1)).
\]
Thus, $b(k)\le x(k)$ for almost all natural numbers $k\in e:=\set{h(n+1)}{n\in d}$.

Let $y\in Y$. Since $\Psi[Y]\les b$,
for almost all natural numbers $k\in e$ we have
\[
\Psi(y)(k)\le b(k)\le x(k),
\]
and thus
\[
\Phi(x,y)(k)\le\Psi(y)(k) \le b(k).
\]
Hence, the set $\Phi[\tilde{X}\x Y]$ is $\les$-bounded on an infinite set,
and thus~\cite[Fact~3.4]{vD} $\les$-bounded.
\epf

As  $\smallcard{(X\sm \tilde{X})\cup\Fin}<\fb$ and the space $Y$ is Hurewicz,  the image $\Phi[((X\sm \tilde{X})\cup\Fin)\x Y]$ is a union of less than $\fb$ sets that are
$\les$-bounded, and is thus $\les$-bounded. Thus, the entire image
$\Phi[(X\cup\Fin)\x Y]$ is $\les$-bounded.
\epf

\section{Cofinal $S$-scales}

For a semifilter $S$, the following special type of $S$-scale
is a natural generalization of the earlier notion of 
cofinal $\bof(S)$-scale~\cite[Definition~2.22]{sfh}.

\bdfn
Let $S$ be a semifilter. A set $X\sub\roth$ with $\card{X}\ge \bof(S)$
is a \emph{cofinal $S$-scale} if for each function $b\in\roth$, we have
\[
b\leS x
\]
for all but less than $\bof(S)$ functions $x\in X$.
\edfn

For example, a set $X\sub\roth$ is a cofinal $\roth$-scale if and only if the set
$X$ is $\fd$-unbounded. Thus, for some semifilters $S$,
cofinal $S$-scales provably exist. But this is not always the case.

\bprp\label{prp:basics}
Let $S$ be a semifilter. 
\be
\item Cofinal $S$-scales are $\le_{S^+}$-unbounded.
\item Every subset, of cardinality $\bof(S)$, of a cofinal $S$-scale is a cofinal $S$-scale.
\item
\label{bd} 
If there is a cofinal $S$-scale, then $\bof(S^+)\le \bof(S)$.
\item
 If $\bof(S)=\fd$, then there is a cofinal $S$-scale~\cite[Lemma~2.23]{sfh}.\qed
\ee
\eprp

\bcor\label{cor:b=d}
Let $F$ be a filter. The following assertions are equivalent:
\be
\item There is a cofinal $F$-scale.
\item $\bof(F)=\bof(F^+)$.
\ee
\ecor
\bpf
$(1)\Impl(2)$: Since $F$ is a filter, we have $F\sub F^+$, and thus $\bof(F)\le\bof(F^+)$.
Apply Proposition~\ref{prp:basics}\eqref{bd}.

$(2)\Impl (1)$: A $\le_{F^+}$-unbounded set $\sset{b_\alpha}{\alpha<\bof(F)}\sub\roth$ is 
$\le_F$-cofinal.
For each ordinal number $\alpha<\bof(F)$, let $x_\alpha\in\roth$ be such that
$\sset{b_\beta,x_\beta}{\beta<\alpha}\le_F x_\alpha$. As $F$ is a filter, the relation
$\le_F$ is transitive, and thus the set
$\sset{x_\alpha}{\alpha<\bof(F)}$ is a cofinal $F$-scale.
\epf

In particular, since $\bof(\cFin)=\fb$ and $\bof(\cFin^+)=\bof(\roth)=\fd$, 
there are cofinal $\cFin$-scales if and only if $\fb=\fd$.

The proof of the following theorem is 
similar to that of Theorems~\ref{thm:master2}--\ref{thm:UU}(1).

\bthm\label{thm:master3}
Let $F$ be a filter and $X\sub\roth$ be a cofinal $F$-scale.
Then, in the realm of hereditarily Lindel\"of spaces,
the real set $X\cup\Fin$ is productively $F$-Menger.\qed
\ethm

\bcor\label{cor:scaleset}
Let $X\sub\roth$ be a cofinal $\cFin$-scale. 
Then the real set $X\cup\Fin$ is productively $F$-Menger for all filters $F$.
\ecor
\bpf
Since the existence of a cofinal $\cFin$-scale implies $\fb=\fd$,
a cofinal $\cFin$-scale is a cofinal $F$-scale for all filters $F$. Apply Theorem~\ref{thm:master3}.
\epf

A \emph{superfilter} (also called \emph{grille} or \emph{coideal})
is a semifilter $S$ such that $a\cup b\in S$ implies $a\in S$ or $b\in S$.
A semifilter $S$ is a superfilter if and only if the semifilter $S^+$ is a filter.
Equivalently, a superfilter is a union of a family of ultrafilters.
For example, the set $\roth=\cFin^+$ is a superfilter.

\bprp\label{prp:cSsEqSs}
Let $S$ be a superfilter. A set $X\sub\roth$ is a cofinal $S$-scale if and only if it is
an $S$-scale.
\eprp
\bpf[Proof $(\Leftarrow)$]
Let $F:=S^+$.
If $b\le_{S^+} c\le_S x$, then $b\le_F c\le_{F^+} x$, and since $F$ is a filter,
we have $b\le_{F^+} x$, that is, $b\leS x$.
\epf

The proof of the following assertion is similar to that of
Proposition~\ref{prp:SscaleSetBdd}.

\bprp\label{prp:NoMoreLabels}
Let $S$ be a semifilter. For each cofinal $S$-scale $X$, the real set $X\cup\Fin$ is $S$-Menger.\qed
\eprp

Let $U$ be an ultrafilter, and $X\sub\roth$ be a $U$-scale. By Proposition~\ref{prp:cSsEqSs},
the set $X$ is in fact a cofinal $U$-scale. Using Proposition~\ref{prp:NoMoreLabels}, we obtain
an alternative derivation of Corollary~\ref{cor:simple}. 
Similarly, Theorem~\ref{thm:master3} generalizes Theorem~\ref{thm:UU}(2).

Theorem~\ref{thm:master3} cannot be extended to all semifilters, and not even to all superfilters:
By Theorems~\ref{thm:bidiNotpM}--\ref{thm:ExistsbiDunb},
the hypothesis $\fd\le\fr$ implies that Theorem~\ref{thm:master3} does not hold for the superfilter $\roth$.
The latter assertion also follows from the following theorem that is, in fact, established by
the proof of Theorem~\ref{thm:strongbidiscale}.

\bthm\label{thm:vsbds}
Assume that $\fb=\fd$. Let $S$ be a semifilter. The following assertions are equivalent:
\be
\item The semifilter $S$ is nonmeager.
\item There are an $S$-scale $X\sub\roth$ and a
$\fd$-concentrated real set $Y$ such that the planar set $(X\cup\Fin)\x Y$ is not Menger.
\item There are a cofinal $S$-scale $X\sub\roth$ and a
$\fd$-concentrated real set $Y$ such that the planar set $(X\cup\Fin)\x Y$ is not Menger.\qed
\ee
\ethm

A related result of Repov\v{s} and Zdomskyy~\cite[Theorem~3.3]{RZ}
asserts that, if $\fb=\fd$, then
there are ultrafilters $U_1$ and $U_2$, a (cofinal) $U_1$-scale $X_1$,
and a (cofinal) $U_2$-scale $X_2$, such that the planar set 
$(X_1\cup\Fin)\x (X_2\cup\Fin)$ is not Menger.

\bthm
Assume that $\fb=\fd$. For every nonmeager filter $F$:
\be
\item In the realm of hereditarily Lindel\"of spaces, there is a 
productively $F$-Menger space that is not Hurewicz and not productively Menger.
\item The property $F$-Menger is strictly between Hurewicz and Menger.
\ee
\ethm
\bpf[Proof of (1)]
By Theorem~\ref{thm:vsbds}(3) and Theorem~\ref{thm:master3},
using that products of Hurewicz sets and $\fd$-concentrated real sets are Menger~\cite[Theorem~4.6]{AddGN}.
\epf

\section{Comments and open problems}

We restrict attention to real sets throughout this section, except for the last subsection.

The \emph{Menger productivity problem}, whether Menger's property is consistently preserved
by finite products, remains open. The hypothesis $\fd\le\fr$ provides two Menger sets
whose product is not Menger (Theorems~\ref{thm:master} and~\ref{thm:bidiNotpM}).
It is well known that this immediately provides a Menger set whose square is not Menger.
Indeed, assume that $X$ and $Y$ are Menger real
sets such that the planar set $X\x Y$ is not Menger. The set $X\cup Y$ is Menger.
We may assume that $X\sub (0,1)$ and $Y\sub (2,3)$. Then the set $X\x Y$
is a closed subset of the square $(X\cup Y)^2$. Menger's property is hereditary for closed subsets.

\subsection{A combinatorial characterization of the cardinal number $\min\{\fr,\fd\}$}
Aubrey~\cite{Aubrey} proved that $\min\{\fd,\fu\}\le\fr$. Since $\fr\le\fu$,
the hypothesis $\fd\le\fr$ in Theorem~\ref{thm:ExistsbiDunb} is equivalent to
the hypothesis $\fd\le\fu$.

Initially, we proved Theorem~\ref{thm:ExistsbiDunb} using a stronger hypothesis.

\bdfn
Let $\bidi$ be the minimal cardinality of a set $X\sub\roth$
such that there is no set $b\in\ici$ with $X\lei b,b\comp$.
\edfn

We observed that $\max\{\fb,\cov(\cM)\}\le\bidi\le\min\{\fr,\fd\}$, and needed that $\bidi=\fd$ to carry out our construction.
It is immediate that $\bidi\le\fd$, and
the argument in the proof of the implication $(2)\Impl(1)$ in Theorem~\ref{thm:ExistsbiDunb} shows that $\bidi\le\fr$.
Answering our question, Mej\'\i{}a pointed out to us that,
by a result of Kamburelis and W\k{e}glorz, 
our upper bound on the cardinal number $\bidi$ is tight~\cite{Mejia} (see Lemma~\ref{lem:mej}).
We thus have the following characterization of $\min\{\fr,\fd\}$.

\bprp\label{prp:bidiBounds}
$\bidi=\min\{\fr,\fd\}$.\qed
\eprp

\subsection{Which $\kappa$-unbounded sets are not productively Menger?}
There are, in ZFC, $\fb$-unbounded sets (e.g., by Proposition~\ref{prp:ExsistScale}).
Every union of less than $\max\{\fb,\fg\}$ Menger sets is Menger~\cite{SF1, sft}.
Since the hypothesis $\fb<\fg$ is consistent, Theorem~\ref{thm:master} and
Corollary~\ref{cor:LnonpM} do not extend to the case $\kappa=\fb$, or
to any cardinal number that is consistently smaller than $\max\{\fb,\fg\}$.

For a $\kappa$-unbounded set, which we may assume to have cardinality $\kappa$,
the present proof of Theorem~\ref{thm:master} requires a partition
$\fd=\Un_{\alpha<\kappa}I_\alpha$ such that for each set $J\sub\kappa$ with $\card{J}<\fd$,
we have $\card{\Un_{\alpha\in J}I_\alpha}<\fd$. It is not difficult to see that
this implies that $\kappa\in\{\cf(\fd),\fd\}$.

\bprb\label{prb:inter}
Assume that $\kappa$ is a cardinal number with $\cf(\fd)<\kappa<\fd$.
Let $X$ be a $\kappa$-unbounded set in $\roth$. Is there necessarily a
$\fd$-concentrated
set $Y$ such that the planar set $X\x Y$ is not Menger?
\eprb

\subsection{Products of Hurewicz sets}
The following problem is intriguing.

\bprb[{Scheepers~\cite[Problem~6.7]{OPiT}}]\label{prb:pH}
Is it consistent that every product of two Hure\-wicz sets is Hurewicz?
\eprb

To this end, it is important to find mild hypotheses implying that there are two
Hurewicz sets whose product is not Hurewicz.
In the notation of Section~\ref{sec:FM}, Menger's property is $\roth$-Menger,
and Hurewicz's property is $\cFin$-Menger.
By Theorems~\ref{thm:bidiNotpM} and~\ref{thm:ExistsbiDunb},
$\roth$-scales need not be productively $\roth$-Menger.
In contrast, by Theorem~\ref{thm:HH},
$\cFin$-scales are productively $\cFin$-Menger.
Thus, modern constructions of Hurewicz sets do not help in regard to
Problem~\ref{prb:pH}.
The classic example of a Hurewicz set (if not counting $\sigma$-compact sets,
which are productively Hurewicz) is a Sierpi\'nski set~\cite{MHP}.
If $\fb=\cov(\cN)=\cof(\cN)$, then there is a $\fb$-Sierpi\'nski set (which is Hurewicz)
whose square is not Hurewicz~\cite[Theorem~43]{cbc}.
No more general constructions are known.

\bprb
Does \CH{} imply that no Sierpi\'nski set is productively Hurewicz?
Productively Menger?
\eprb

By Theorem~\ref{thm:pMispH}, if $\fb=\fd$, then every productively Menger set
is productively Hurewicz. By the discussion following that theorem, if $\fb<\fg$
then there are productively Menger sets that are not Hurewicz.

\bprb
Assume \CH{}.
Do the classes of productively Menger and productively Hurewicz sets coincide?
\eprb

Recall that for meager semifilters $S$, being $S$-Menger is equivalent to being
Hurewicz~\cite[Theorem~2.32]{sfh}. By Theorem~\ref{thm:HH},
in this case, for each $S$-scale $X$, the set $X\cup\Fin$ is productively $S$-Menger.

\bprb\label{prb:jj}
Assume \CH{}.
For which semifilters $S$ there is 
an $S$-scale $X$ such that the set $X\cup\Fin$ is not productively $S$-Menger?
\eprb

Thus meager semifilters do not have the property in Problem~\ref{prb:jj}.
By Theorem~\ref{thm:master3} and Proposition~\ref{prp:cSsEqSs}, 
ultrafilters are also not in that category. But the full semifilter $\roth$ 
is in this category, by Theorem~\ref{thm:vsbds}.

A $\fb$-scale~\cite{MHP}
is a particularly simple kind of a $\cFin$-scale, for $\cFin$ the filter of cofinite sets.

\bprb
Let $X\sub\roth$ be a $\fb$-scale. Is the set $X\cup\Fin$ necessarily productively Menger?
\eprb

If $\fu<\fg$, then every $\fd$-concentrated set (in particular, every
union of an $S$-scale, for some semifilter $S$, and the set $\Fin$) 
is productively Menger~\cite[Theorem~4.7]{SPMProd}.

\subsection{Strictly unbounded sets}
Say that a set $X\sub\roth$ is \emph{strictly unbounded}
if for every set $A\sub\roth$ of cardinality smaller
than $\fd$ there is a function $x\in X$ such that
$A\lei x$.
Let $X\sub\roth$ be a strictly unbounded set.
By the argument in the proof of Theorem~\ref{thm:pMispH}(1),
the set $X$ contains a $\fd$-unbounded set.
By Theorem~\ref{thm:master}, there is a
$\fd$-concentrated set $Y$
such that the planar set $X\x Y$ is not Menger.
If $\fb=\fd$, then every unbounded set in $\roth$ is strictly unbounded.
We thus obtain a generalization of Theorem~\ref{thm:pMispH}.

The construction in Theorem~\ref{thm:ExistsbiDunb} that provides a Menger
set that is not productively Menger provides, in fact, a Menger
strictly unbounded set.

A negative answer to the second item of the following problem implies
a negative solution for the Menger productivity problem.

\bprb
\mbox{}
\be
\item Is it consistent that $\fr<\fd$ and there are strictly unbounded Menger sets?
\item Is it consistent that there are no strictly unbounded Menger sets?
\ee
\eprb

\subsection{General spaces}
Let $S$ be a semifilter.
Restricting the definition of \emph{$S$-Menger} spaces to countable open covers, we obtain the definition of \emph{countably $S$-Menger} spaces. 
This makes it possible to eliminate
the adjective ``hereditarily Lindel\"of'' in most of our theorems.

For \emph{general} Hurewicz spaces, the following problem remains open, even for
the so-called $\fb$-scales~\cite[Definition~2.8]{MHP}.

\bprb
Let $\cFin$ be the semifilter of cofinite sets,
$X\sub\roth$ be an $\cFin$-scale,
and $Y$ be a Hurewicz space.
Is the product space $(X\cup\Fin)\x Y$ necessarily Hurewicz?
\eprb

Theorem~\ref{thm:HH} provides a positive answer for hereditarily Lindel\"of spaces $Y$. 
But this restriction is only needed for deducing that the product space
$(X\cup\Fin)\x Y$ is Lindel\"of (and similarly for the other results in Section~\ref{sec:FM}).
A positive solution for the following problem would suffice.

\bprb
Let $X$ be a real set of cardinality smaller than $\fb$, 
and $Y$ be a Hurewicz space.
Is the product space $X\x Y$ necessarily Lindel\"of?
\eprb

\ed

%% file: shortcuts.tex
\newcommand{\nc}{\newcommand}
\DeclareMathOperator{\PS}{P}
\nc{\thusfar}{\par\bigskip\centerline{\my{--- Edited thus far ---}}\par\bigskip}
\nc{\lei}{\le^\oo}
\nc{\card}[1]{\left|#1\right|}
\nc{\medcard}[1]{\biggl|\,#1\,\biggr|}
\nc{\smallcard}[1]{|\,#1\,|}
\nc{\bds}{bidirectional $\roth$-scale}
\nc{\bbN}{\mathbb{N}}
\nc{\beq}{\begin{eqnarray*}}\nc{\eeq}{\end{eqnarray*}}
\nc{\mbq}{\mb{?}}
\nc{\mb}[1]{{\mbox{\textbf{#1}}}}
\nc{\nop}{$\times$}
\nc{\fbn}{\!\!\fbox{\!\nop\!}\!\!}
\nc{\yup}{\checkmark}
\nc{\forces}{\Vdash}
\nc{\name}[1]{\dot{#1}}
\nc{\tf}{\my{FINISHED THUS FAR}}
\nc{\FU}{Fr\'echet--Urysohn}
\nc{\gs}{$\gamma$~space}
\nc{\Ga}{\Gamma}\nc{\Om}{\Omega}
\nc{\smallbinom}[2]{\begin{psmallmatrix} #1\\ #2 \end{psmallmatrix}}
\nc{\bgamma}{\smallbinom{\Om}{\Ga}}
\newcommand{\two}{\{0,1\}}
\nc{\productive}[2]{\bigl(#1,\allowbreak #2\bigr)^\x}
\nc{\Sel}{\mathsf{S}}
\nc{\sset}[2]{\{\,#1 : #2\,\}}
\nc{\smb}[1]{{\!\!\mb{#1}\!\!}}
\nc{\medset}[2]{{\biggl\{\,#1 : #2\,\biggr\}}}
\nc{\smallmedset}[2]{{\bigl\{\,#1 : #2\,\bigr\}}}
\nc{\set}[2]{{\left\{\,#1 : #2\,\right\}}}
\nc{\seq}[2]{{\la\, #1 : #2\,\ra}}
\nc{\eseq}[1]{#1_1, \allowbreak #1_2, \allowbreak\dots} 
\nc{\cube}{(\Cantor)^\bbN}
\nc{\Match}{\op{Match}}
\nc{\concat}[1]{\hat{\phantom{a}}\langle #1\rangle}
\nc{\poset}{\mathbb{P}}
\nc{\fn}[1]{{\op{Fn}(#1\times\w,2)}}
\nc{\linadd}{\op{linadd}}
\nc{\nonprod}{\non^\x}
\nc{\alephes}{{\aleph_0}}
\nc{\my}[1]{\marginpar{\textcolor{red}{***}}\textcolor{red}{#1}}
\nc{\later}[1]{{\color{green} #1}}
\nc{\BTs}[1]{{\color{green} #1 (BT)}}
\nc{\Cp}{\op{C}_\mathrm{p}}
\nc{\Bp}{\op{B}_p}
\nc{\Pa}[8]{\bibitem{#1} {#2}, \emph{#3}, {#4} \textbf{#5} ({#6}), {#7}--{#8}.}
\nc{\tPa}[5]{\bibitem{#1} {#2}, \emph{#3}, {#4}, to appear.}
\nc{\sPa}[4]{\bibitem{#1} {#2}, \emph{#3}, {#4}, submitted.}
\nc{\Bc}[9]{\bibitem{#1} {#2}, \emph{#3}, in: \textbf{#4} (#5), #6 #7, #8--#9.}
\nc{\fD}{\mathfrak{D}}
\nc{\fX}{\mathfrak{X}}
\nc{\Onbd}{\Op_{\mathrm{nbd}}} 
\nc{\Omnb}{\Om_{\mathrm{nbd}}} 
\nc{\od}{\mathfrak{od}}
\nc{\Setting}[7]{\xymatrix@R=4pt@C=7pt{#1\ar@{-}[r]&#2\ar@{-}[r]&#3\\&#4\ar@{-}[u]\\
#5\ar@{-}[uu]\ar@{-}[r] & #6\ar@{-}[u]\ar@{-}[r] & #7\ar@{-}[uu]}}
\nc{\mx}[1]{\begin{matrix}#1\end{matrix}}
\nc{\plim}{p\txt{-}\lim}
\nc{\Bgp}{{\Z^\bbN}}
\nc{\Cgp}{{{\Z_2}^\bbN}}
\nc{\Cite}[1]{\textbf{[#1]}}
\nc{\Next}[1]{{#1^+}}
\nc{\cFin}{\mathrm{cF}}
\nc{\intvl}[2]{{[#1(#2),\allowbreak #1(#2\!+\!1))}}
\nc{\Bdd}{\mathbf{B}}
\nc{\Dfin}{\mathfrak{D}_\mathrm{fin}}
\nc{\grbl}{{\mbox{\textit{\tiny gp}}}}
\nc{\bbP}{\mathbb{P}}
\nc{\BOfat}{\B_{\Om_{\mathrm{fat}}}}
\nc{\Bgood}{\B_{\mathrm{good}}}
\nc{\compactN}{\cl{\mathbb{N}}}
\nc{\blocks}[2]{\op{cl}_{#2}(#1)}
\nc{\blocksplus}[2]{\op{cl}^+_{#2}(#1)}
\nc{\arx}[1]{\texttt{http://arxiv.org/math/#1}}
\nc{\bq}{\begin{quote}}
\nc{\eq}{\end{quote}}
\nc{\cl}[1]{\overline{#1}}
\nc{\CH}{the Continuum Hypothesis}
\nc{\MA}{Martin's Axiom}
\nc{\Bfat}{\B_\mathrm{fat}}
\nc{\inv}{^{-1}}
\nc{\Cantor}{{\two^\bbN}}
\nc{\bP}{\mathbf{P}}
\nc{\bof}{\op{\fb}}
\nc{\dof}{\op{\fd}}
\nc{\bofF}{\bof(\cF)}
\nc{\sr}[3]{\underset{\mbox{#3}}{\mbox{#1}}}
\nc{\gp}{\binom{\Om}{\Ga}}
\nc{\gpsmall}{\mbox{$\gp$}}
\nc{\gig}{\gimel}
\nc{\gns}{\sone(\Om,\gig)}
\nc{\nsr}[2]{#1}
\nc{\Srg}{{\mathbb{S}}}
\nc{\Srgs}{{\mathbb{S}^*}}
\nc{\NN}{{\bbN^{\bbN}}}
\nc{\ZN}{{\Z^{\bbN}}}
\nc{\NNup}{{\bbN^{\uparrow\bbN}}}
\nc{\Pof}{\op{P}}
\nc{\PN}{{\Pof(\bbN)}}
\nc{\roth}{{[\bbN]^{\mbox{\tiny $\infty$}}}} 
\nc{\Fin}{[\bbN]^{\text{$<\!\!\infty$}}} 
\nc{\ici}{[\bbN]^{ \infty, \infty}}
\nc{\Inc}{{\compactN^{\uparrow\bbN}}}
\nc{\powInc}[1]{{\big(\Inc\big)^{#1}}}
\nc{\powFin}[1]{{\big(\Fin\big)^{#1}}}
\nc{\powPN}[1]{{\big(\PN\big)^{#1}}}
\nc{\NcompactN}{{\compactN^\bbN}}
\nc{\Uarrow}{\smash{\big\uparrow}}
\nc{\LE}{\preccurlyeq}
\nc{\GE}{\succcurlyeq}
\nc{\op}{\operatorname}
\nc{\im}{\op{im}}
\nc{\Span}{\op{span}}
\nc{\maxfin}{\op{maxfin}}
\nc{\ran}{\op{range}}
\nc{\iso}{\cong}
\nc{\Madd}{{\M}^\star}
\nc{\cI}{\mathcal{I}}
\nc{\cJ}{\mathcal{J}}
\nc{\scrA}{\mathscr{A}}
\nc{\scrB}{\mathscr{B}}
\nc{\scrC}{\mathscr{C}}
\nc{\scrD}{\mathscr{D}}
\nc{\scrF}{\mathscr{F}}
\nc{\scrK}{\mathscr{K}}
\nc{\A}{\forall}
\nc{\B}{\mathrm{B}}
\nc{\cB}{\mathcal{B}}
\nc{\bB}{\mathbf{B}}
\nc{\BS}{\mathbf{B}(\mathcal{S})}
\nc{\BF}{\mathbf{B}(\mathcal{F})}
\nc{\BU}{\mathbf{B}(\mathcal{U})}
\nc{\cSp}{\mathcal{S}^+}
\nc{\cFp}{\mathcal{F}^+}
\nc{\cUp}{\mathcal{U}^+}

\nc{\BG}{\B_\Ga}
\nc{\BL}{\B_\Lambda}
\nc{\BT}{\B_\Tau}
\nc{\BTstar}{\B_{\Tau^*}}
\nc{\BO}{\B_\Om}
\nc{\DO}{\cD_\Om}
\nc{\KO}{\cK_\Om}
\nc{\CG}{C_\Ga}
\nc{\CL}{C_\Lambda}
\nc{\CT}{C_\Tau}
\nc{\CTstar}{C_{\Tau^*}}
\nc{\CO}{C_\Om}
\nc{\COgp}{C_{\Om^{\grbl}}}
\nc{\CLgp}{C_{\Lambda^{\grbl}}}
\nc{\BOgp}{\B_{\Om}^{\grbl}}
\nc{\BLgp}{\B_{\Lambda^{\grbl}}}
\nc{\sfC}{\mathsf{C}}
\nc{\sfD}{\mathsf{D}}
\nc{\bD}{\mathbf{D}}
\nc{\Tau}{\mathrm{T}}
\nc{\cA}{\mathcal{A}}
\nc{\cK}{\mathcal{K}}
\nc{\cD}{\mathcal{D}}
\nc{\cF}{\mathcal{F}}
\nc{\cS}{\mathcal{S}}
\nc{\cT}{\mathcal{T}}
\nc{\cG}{\mathcal{G}}
\nc{\cY}{\mathcal{Y}}
\nc{\J}{\mathcal{J}}
\nc{\cL}{\mathcal{L}}
\nc{\cM}{\mathcal{M}}
\nc{\cN}{\mathcal{N}}
\nc{\cH}{\mathcal{H}}
\nc{\cO}{\mathcal{O}}
\nc{\Op}{\mathrm{O}}
\nc{\rmA}{\mathrm{A}}
\nc{\rmF}{\mathrm{F}}
\nc{\rmB}{\mathrm{B}}
\nc{\rmD}{\mathrm{D}}
\nc{\rmP}{\mathrm{P}}
\nc{\cC}{\mathcal{C}}
\nc{\cP}{\mathcal{P}}
\nc{\bbQ}{\mathbb{Q}}
\nc{\bbR}{\mathbb{R}}
\nc{\cU}{\mathcal{U}}
\nc{\Un}{\bigcup}
\nc{\cV}{\mathcal{V}}
\nc{\cW}{\mathcal{W}}
\nc{\Z}{{\mathbb Z}}
\nc{\Impl}{\Rightarrow}
\long\def\forget#1\forgotten{\marginpar{\textcolor{green}{Forgetting...}}}
\nc{\ft}{\mathfrak{t}}
\nc{\fb}{\mathfrak{b}}
\nc{\fc}{\mathfrak{c}}
\nc{\fd}{\mathfrak{d}}
\nc{\fg}{\mathfrak{g}}
\nc{\oo}{\infty}
\nc{\fr}{\mathfrak{r}}
\nc{\fk}{\mathfrak{k}}
\nc{\bidi}{\mathfrak{bidi}}
\nc{\fu}{\mathfrak{u}}
\nc{\fh}{\mathfrak{h}}
\nc{\fp}{\mathfrak{p}}
\nc{\fj}{\mathfrak{j}}
\nc{\fs}{\mathfrak{s}}
\nc{\w}{\omega}
\nc{\x}{\times}
\nc{\Iff}{\Leftrightarrow}
\newcommand\comp{^{\text{\tt c}}}
\nc{\nin}{\notin}
\nc{\cat}{\hat{\ }}
\nc{\sub}{\subseteq}
\nc{\spst}{\supseteq}
\nc{\sm}{\setminus}
\nc{\as}{\subseteq^*}
\nc{\les}{\le^*}
\nc{\leinf}{\le^{\infty}}
\nc{\leS}{\le_S}
\nc{\leF}{\le_{\mathcal{F}}}
\nc{\leU}{\le_{\mathcal{U}}}
\nc{\rest}{\restriction}

\nc{\la}{\langle}
\nc{\ra}{\rangle}
\nc{\E}{\exists}
\nc{\dom}{\op{dom}}
\nc{\cov}{\op{cov}}
\nc{\add}{\op{add}}
\nc{\cof}{\op{cof}}
\nc{\cf}{\op{cf}}
\nc{\non}{\op{non}}
\nc{\unif}{\op{non}}
\nc{\COV}{\op{COV}}
\nc{\ADD}{\op{ADD}}
\nc{\COF}{\op{COF}}
\nc{\NON}{\op{NON}}
\nc{\impl}{\to}
\nc{\Lp}{\mathcal{L_\p}}
\nc{\Wlog}{without loss of generality}
\newtheorem{thm}{Theorem}[section]
\nc{\bthm}{\begin{thm}} \nc{\ethm}{\end{thm}}
\newtheorem{prop}[thm]{Proposition}
\nc{\bprp}{\begin{prop}} \nc{\eprp}{\end{prop}}
\newtheorem{fact}[thm]{Fact}
\nc{\bfct}{\begin{fact}} \nc{\efct}{\end{fact}}
\newtheorem{prob}[thm]{Problem}
\nc{\bprb}{\begin{prob}} \nc{\eprb}{\end{prob}}
\newtheorem{lem}[thm]{Lemma}
\nc{\blem}{\begin{lem}} \nc{\elem}{\end{lem}}
\newtheorem{claim}[thm]{Claim}
\nc{\bclm}{\begin{claim}} \nc{\eclm}{\end{claim}}
\newtheorem{cor}[thm]{Corollary}
\nc{\bcor}{\begin{cor}} \nc{\ecor}{\end{cor}}
\newtheorem{conj}[thm]{Conjecture}
\nc{\bcnj}{\begin{conj}} \nc{\ecnj}{\end{conj}}
\theoremstyle{definition}
\newtheorem{defn}[thm]{Definition}
\nc{\bdfn}{\begin{defn}} \nc{\edfn}{\end{defn}}
\newtheorem{obs}[thm]{Observation}
\nc{\bobs}{\begin{obs}} \nc{\eobs}{\end{obs}}
\theoremstyle{remark}
\newtheorem{rem}[thm]{Remark}
\nc{\brem}{\begin{rem}} \nc{\erem}{\end{rem}}
\newtheorem{cnv}[thm]{Convention}
\nc{\bcnv}{\begin{cnv}} \nc{\ecnv}{\end{cnv}}
\newtheorem{exam}[thm]{Example}
\nc{\bexm}{\begin{exam}} \nc{\eexm}{\end{exam}}
\nc{\bpf}{\begin{proof}} \nc{\epf}{\end{proof}}
\nc{\be}{\begin{enumerate}}
\nc{\ee}{\end{enumerate}}
\nc{\bi}{\begin{itemize}}
\nc{\bimy}{\my{\begin{itemize}}
\nc{\eimy}{\end{itemize}}}
\nc{\itm}{\item}
\nc{\ei}{\end{itemize}}
\nc{\Subsection}[1]{\goodbreak\subsection*{#1}}

\nc{\sone}{\mathsf{S}_1}
\nc{\sfin}{\mathsf{S}_\mathrm{fin}}
\nc{\ufin}{\mathsf{U}_\mathrm{fin}}
\nc{\Split}{\mathsf{Split}}

\nc{\gone}{\mathsf{G}_1}    \nc{\gfin}{\mathsf{G}_\mathrm{fin}}

%% file: pMFinal.bbl
\begin{thebibliography}{00}

\bibitem{Aubrey} J. Aubrey,
\emph{Combinatorics for the dominating and unsplitting numbers},
Journal of Symbolic Logic \textbf{69} (2004), 482--498.

\bibitem{AurD} L. Aurichi,
\emph{D-spaces, topological games, and selection principles},
Topology Proceedings \textbf{36} (2010), 107--122.

\bibitem{AurTall} L. Aurichi, F. Tall, 
\emph{Lindel\"of spaces which are indestructible, productive, or D}, 
Topology and its Applications \textbf{159} (2012), 331--340.

\bibitem{BZCSFQdM} T. Banakh, L. Zdomskyy,
\emph{The Coherence of Semifilters: a Survey},
in: \textbf{Selection Principles and Covering Properties in Topology} (L. Ko\v{c}inac, ed.),
Quaderni di Matematica \textbf{18}, Seconda Universit\`a di Napoli, Caserta, 2006, 53--99.

\bibitem{BaSh01}
T. Bartoszy\'nski, S. Shelah,
\emph{Continuous images of sets of reals},
Topology and its Applications \textbf{116} (2001), 243--253.

\bibitem{Ideals} {T. Bartoszy\'nski, B. Tsaban},
\emph{Hereditary topological diagonalizations and the Menger-Hurewicz Conjectures},{Proceedings of the American Mathematical Society} \textbf{134} ({2006}), {605}--{615}.

\bibitem{BlassHBK} A. Blass,
\emph{Combinatorial cardinal characteristics of the continuum},
in: \textbf{Handbook of Set Theory} (M. Foreman, A. Kanamori, eds.),
Springer, 2010, 395--489.

\bibitem{Bri} W. Brian,
Answer to \emph{A property of the Frechet filter and every ultrafilter},
MathOverflow, 2015.
\url{http://mathoverflow.net/questions/201171}

\bibitem{ChoZdo} D. Chodounsky, D. Repov\v{s}, L. Zdomskyy,
\emph{Mathias forcing and combinatorial covering properties of filters},
Journal of Symbolic Logic, Journal of Symbolic Logic \textbf{80} (2015), 1398--1410. 

\bibitem{vD}
E. van Douwen,
\emph{The integers  and  topology},
in: \textbf{Handbook of Set Theoretic Topology} (K. Kunen, J. Vaughan, editors),
North-Holland, Amsterdam, 1984, 111--167.

\bibitem{Hure25} W. Hurewicz,
\emph{\"Uber eine Verallgemeinerung des Borelschen Theorems},
Mathematische Zeitschrift \textbf{24} (1925), 401--421.

\bibitem{coc2} W. Just, A. Miller, M. Scheepers, P. Szeptycki,
\emph{The combinatorics of open covers II},
Topology and its Applications \textbf{73} (1996), 241--266.

\bibitem{KW} A. Kamburelis, B. W\k{e}glorz, 
\emph{Splittings}, Archive for Mathematical Logic \textbf{35} (1996), 263--277.

\bibitem{Ash} A. Kumar, Answer to
\emph{A classic cardinal characteristic of the continuum in disguise?},
MathOverflow, 2015.
\url{http://mathoverflow.net/questions/201170}

\bibitem{Law} L. Lawrence,
\emph{The influence of a small cardinal on the product of a Lindel\"of space and the irrationals},
Proceedings of the American Mathematical Society \textbf{110} (1990), 535--542.

\bibitem{Mejia}
D. Mej\'\i{}a,
Answer for
\emph{Bidi: A new cardinal characteristic of the continuum?},
MathOverflow, 2015.
\url{http://mathoverflow.net/questions/206348}

\bibitem{Menger24} K. Menger,
\emph{Einige \"Uberdeckungss\"atze der Punktmengenlehre},
Sitzungsberichte der Wiener Akademie \textbf{133} (1924), 421--444.

\Pa{SPMProd}{A. Miller, B. Tsaban, L. Zdomskyy}{Selective covering properties of product spaces}{Annals of Pure and Applied Logic}{165}{2014}{1034}{1057}

\bibitem{RZ} D. Repov\v{s}, L. Zdomskyy,
\emph{On M-separability of countable spaces and function spaces}, 
Topology and its Applications \textbf{157} (2010), 2538--2541.

\Pa{cbc}{M. Scheepers, B. Tsaban}{The combinatorics of Borel covers}{Topology and its Applications}{121}{2002}{357}{382}

\bibitem{Talagrand} M. Talagrand,
\emph{Compacts de fonctiones mesurables et filtres non mesurables},
Studia Mathematica \textbf{67} (1980), 13--43.

\bibitem{Tod95} S. Todor\v{c}evi\'c,
\emph{Aronszajn orderings},
Publications de l'Insitut Mathematique \textbf{57} (1995), 29--46.

\Bc{OPiT}{B. Tsaban}{Selection principles and special sets of reals}{Open Problems in Topology II}{E. Pearl, ed.}{Elsevier B.V.}{2007}{91}{108}

\bibitem{MHP} B. Tsaban,
\emph{Menger's and Hurewicz's Problems: Solutions from ``The Book'' and refinements},
in: Set Theory and its Applications (L. Babinkostova, A. Caicedo, S. Geschke, M. Scheepers, editors),
Contemporary Mathematics \textbf{533} (2011), 211--226.

\bibitem{AlgSelRT} B. Tsaban,
\emph{Algebra in the Stone--\v{C}ech compactification, selections, and additive Ramsey theory},
submitted for publication.

\Pa{sft}{B. Tsaban, L. Zdomskyy}{Combinatorial images of sets of reals and semifilter trichotomy}{Journal of Symbolic Logic}{73}{2008}{1278}{1288}

\Pa{sfh}{B. Tsaban, L. Zdomskyy}{Scales, fields, and a problem of Hurewicz}{Journal of the European Mathematical Society}{10}{2008}{837}{866}

\Pa{AddGN}{B. Tsaban, L. Zdomskyy}{Additivity of the Gerlits--Nagy property and concentrated sets}{Proceedings of the American Mathematical Society}{142}{2014}{2881}{2890}

\bibitem{SF1} L. Zdomskyy,
\emph{A semifilter approach to selection principles},
Commentationes Mathematicae Universitatis Carolinae \textbf{46} (2005), 525--539.

\end{thebibliography}
